\documentclass[times,sort&compress,3p]{elsarticle}
\journal{Journal of Multivariate Analysis}
\usepackage{amsmath, amssymb, enumerate, amsthm,color}

\def\esp{{\rm E} }
\def\R{\mathbb{R}}
\def\tr{\,\textmd{tr}\,}
\def\<{\langle}
\def\>{\rangle}
\def\rank{\mathrm{rank}\, }
\begin{document}

\begin{frontmatter}
\title{The Laplace transform $(\det s)^{-p}\exp \mathrm{tr}( s^{-1}w)$ and \\
the existence of non-central Wishart distributions}

\author[GL]{G\'erard Letac}
\ead{gerard.letac@math.univ-toulouse.fr} 
\author[HM]{H\'el\`ene Massam}
\ead{massamh@yorku.ca}
\address[GL]{Institut de math\'ematiques de Toulouse, Universit\'e Paul-Sabatier, 31062, Toulouse,
France}
\address[HM]{Department of Mathematics and Statistics, York
University, Toronto, Ontario, Canada M3J 1P3}  

 \begin{abstract} 
The problem considered in this paper is to find when  the non-central Wishart distribution, defined on the cone $\overline{\mathcal{P}_d}$ of  positive semidefinite matrices of order $d$ and with a  real-valued  shape parameter $p$, does   exist.
This can be  reduced  to the study  of the measures $m(n,k,d)$ defined on $\overline{\mathcal{P}_d}$ and with Laplace transform $(\det s)^{-n/2}\exp \tr(s^{-1}w)$, where $n$ is an integer and $w=\mathrm{diag}(0,\ldots,0,1,\ldots,1)$ has order $d$ and rank $k.$  Our two main results are the following: we compute $m(d-1,d,d)$ and we show  that neither $m(d-2,d,d)$ nor  $m(d-2,d-1,d)$ exists. These  facts  solve the problems  of the  existence and computation of these non-central Wishart distributions.
\end{abstract}


%\noindent AMS classification: primary 60H10; secondary 60B11. 

%\bigskip \noindent \textsc{Abbreviated title:} Non-central Wishart laws. 

\end{frontmatter}

\section{Introduction} The non-central Wishart distribution is traditionally defined as the distribution of  the random symmetric real matrix $X=Y_1Y_1^\top+\cdots+Y_nY_n^\top$, where $Y_1 \ldots, Y_n \in \mathbb{R}^d$ are independent Gaussian column vectors  with the same non-singular  covariance matrix $\Sigma$ and  respective means $m_1, \ldots, m_n$ which  are not necessarily equal.  For $s$  in the open cone $\mathcal{P}_d$ of positive definite symmetric matrices of order $d$ and 
$w=m_1m_1^\top+\cdots+m_nm_n^\top$  in the closed cone of  positive semidefinite matrices $\overline{\mathcal{P}_d}$, one can readily derive the Laplace transform
\begin{equation}
\label{NCWstat1}{\rm E} (e^{-\tr (sX)})=\frac{1}{\det (I_d+2\Sigma s)^{n/2}} \, e^{-\tr \{2s(I_d+2\Sigma	s)^{-1}w\} }. 
\end{equation}
It is important to note that in this formula, the rank $k$ of $w$ is such that $k \le \min (n, d)$.

In the same way that the familiar chi-square distribution with $n$ degrees of freedom can be extended to the gamma distribution with a continuous shape parameter $p$ by replacing the half integer $n/2$ by $p\in (0, \infty)$, it is tempting to extend the values that the power of $\det (I_d+2\Sigma	s)$ can take in \eqref{NCWstat1}. The question is then: given $\Sigma\in \mathcal{P}_d$ and $w\in \overline{\mathcal{P}_d}$, for which values of  $p>0$ does there exist a probability distribution on $\overline{\mathcal{P}_d}$ for $X$ such that, for all $s\in \mathcal{P}_d$, we have 
\begin{equation}
\label{NCW}
{\rm E} (e^{-\tr (sX)})=\frac{1}{\det (I_d+2\Sigma s)^{p}} \, e^{-\tr \{ 2s(I_d+2\Sigma	s)^{-1}w \}}. 
\end{equation}
The hypothetical distribution for $X$ satisfying \eqref{NCW} can be called a non-central Wishart distribution with parameters $(2p,w,\Sigma),$ or $\mathcal{NCW}(2p,w,\Sigma)$
for short.  

The first question is to determine what are the possible values of $p$. For $w=0$, i.e., for the extension of the ordinary Wishart distributions, the problem is solved by Gindikin's theorem \cite{Gindikin}, which says that  $\mathcal{NCW}(2p,0,\Sigma)$ exists if and only if
\begin{equation}
\label{GYN}
p \in \Lambda_d = \left\{\frac{1}{2}\, ,\frac{2}{2}\, ,\ldots,\frac{d-2}{2} \right\} \cup \left[(d-1)/2,\infty \right).
\end{equation}
 
In Proposition~8 below, we give  a  proof of Gindikin's Theorem, which we borrow from \cite{Peddada}  because it is the most adapted to the zonal polynomial techniques used in this paper. The Wishart distributions with $w=0$ and when  $2p$ not necessarily an integer and  satisfying $2p>d-1$ have densities  proportional to $e^{-\tr (\Sigma^{-1}x)}(\det x)^{p-(d+1)/2}$ concentrated on $\mathcal{P}_d.$ A simple  argument (see Proposition~1 below) shows that  if $\mathcal{NCW}(2p,w,\Sigma)$ exists, then for all $\epsilon\geq 0$  the distribution  $\mathcal{NCW}(2p,\epsilon w,\Sigma)$ also exists. By a limiting argument, this observation implies that  if $\mathcal{NCW}(2p,w,\Sigma)$ exists, then  $\mathcal{NCW}(2p,0,\Sigma)$ exists. Hence $p\in \Lambda_d$ and the  difficulty concerning the possible values of $p$ is overcome. 

The second question is a  more difficult  one. In the ordinary non-central Wishart \eqref{NCWstat1}, we have seen that the rank $k$ of $w$ must be such that $k \le n=2p$. Do we have a similar constraint   for  the $\mathcal{NCW}(2p,w,\Sigma)$  distribution?  In other words, 
 does the $\mathcal{NCW}(2p,w,\Sigma)$ defined by  \eqref{NCW} fail to exist if $k > 2p?$ We have a counterexample when $(d-1)/2<p< d/2$ and $k=d.$ Indeed, similarly to the case $w=0$, it has been proved (see, e.g., \cite{LM2008}) that, for any $w$, the distribution 
$\mathcal{NCW}(2p,w,\Sigma)$ exists and has a density if  $(d-1)/2<p$. 

The main objective of the present paper is to answer this question. We will see that: 
\begin{itemize} 
\item [(a)] 
If $p\in \left\{ {1}/{2}, {2}/{2},\ldots, (d-2)/{2} \right\}$, then $\mathcal{NCW}(2p,w,\Sigma)$ exists if and only if the rank $k$ of $w$ is such that $k \leq 2p$. 
\item [(b)]
If $p\geq (d-1)/2$, then  $\mathcal{NCW}(2p,w,\Sigma)$ exists for any $w.$ 
\end{itemize}

As just noted above, Statement~(b) has already been proved for $p>(d-1)/2$ and  the density of  the $\mathcal{NCW}(2p,w,\Sigma)$ distribution is known. For $p=(d-1)/2$, Statement~(b) is established by observing that
$$
  \mathcal{NCW}(d-1,w,\Sigma)=\lim_{p\downarrow (d-1)/2}\mathcal{NCW}(2p,w,\Sigma).
$$
  
\section{Reduction of the problem: The measures $m(2p,k,d)$}
\label{sec:2}

 Let $k$ be an integer such that $0\leq k\leq d$. We consider the  diagonal matrix $I(k,d)$ with its first $d-k$
 diagonal terms  equal to $0$ and the last $k$   equal to 1, viz.
$$
I(k,d)=\left[\begin{array}{cc}0_{d-k}&0\\0&I_k\end{array}\right].
$$
For $p\in \Lambda_d$ we define the positive measure $m(2p,k,d)$  on $\overline{\mathcal{P}_d}$ such that, for all $s\in  \mathcal{P}_d$, we have \begin{equation}
\label{MW}
\int_{\overline{\mathcal{P}_d}}e^{-\tr(sx)}m(2p,k,d)(dx)=\frac{1}{(\det s)^p} \, e^{\tr \{s^{-1} I(k,d)\} }.
\end{equation}

Note that $m(2p,k,d)$ may or may not exist. For example, formula \eqref{m111} below shows that the density of  $m(1,1,1)$ on $(0,\infty)$  is 
${\cosh (2\sqrt{x})}/{\sqrt{\pi x}}$. 
More generally with $ p>0$, we have
$$
m(0,0,1)=\delta_0, \quad m(0,1,1)=\delta_0+\left\{ \sum_{n=1}^{\infty}\frac{x^{n-1}}{n!(n-1)!}\right\} \mathbf{1}_{(0,\infty)}(x)dx,
$$
$$
m(p,0,1)=\frac{x^{p-1}}{\Gamma(p)} \, \mathbf{1}_{(0,\infty)}(x)dx, \quad m(p,1,1)=\left\{ \sum_{n=0}^{\infty}\frac{x^{n+p-1}}{n!\Gamma(n+p)}\right\} \mathbf{1}_{(0,\infty)}(x)dx,
$$
where $\delta_0$ denotes Kronecker's delta.

For $2p>d-1$, formula \eqref{LTM2} below gives  $m(2p,d,d).$ 
If $k$ and $n$ are integers such that $0\leq k\leq n\leq d$, formula \eqref{SS} gives $m(n,k,d)$. 
Finally $m(1,2,2)$ is computed in Section~\ref{sec:3}. The measure $m(d-1,d,d)$  is computed in Section~\ref{sec:4.3}  and details about $m(2,3,3)$  are given in Section~\ref{sec:4.4}. We will show below that these examples are the only cases of existence. For instance the function $s\mapsto \exp \tr(s^{-1})$ on $\mathcal{P}_d$ is not the Laplace transform of a positive measure if $d\geq 2.$ 

The following proposition links this  unbounded measure $m(2p,k,d)$ defined by \eqref{MW} with our initial existence problem. Proposition~1  is important for the solution to the problem of  existence or  non existence of the non-central  Wishart with continuous shape parameter $p$: it focuses the problem at its core by ignoring the normalization constant and  the parameter $\Sigma$, and by  reducing the parameter $w$ to its most important characteristic, namely its rank $k$.

\bigskip \noindent % Formerly 2.1
\textbf{Proposition 1.} 
{\it Let   $\Sigma\in \mathcal{P}_d,$  $w\in \overline{\mathcal{P}_d}$ and $p\in\Lambda_d.$ Suppose that $\rank w=k.$ Then $\mathcal{NCW}(2p,w,\Sigma)$ as defined by \eqref{NCW} exists if and only if $m(2p,k,d)$ exists. Furthermore, if  $m(2p,k,d)$ exists, then $p$ belongs to the Gindikin set \eqref{GYN}.} 

\bigskip 
\noindent 
\textit{Proof.} Assume that $m(2p,k,d)$ exists and let us show that $\mathcal{NCW}(2p,w,\Sigma)$ exists. The proof is based on the following principle. Let $\mu$ be a positive measure on a finite-dimensional real linear space $E$ such that its Laplace transform $L_{\mu}(s)=\int_E e^{-\<s,x\>}\mu(dx)$ is finite on some convex  subset  $D(\mu)$ of the dual space $E^\top$ with a non-empty interior. Let $a$ be a linear automorphism of $E^\top$ and let $b\in E^\top$ such that  $L_{\mu}\{a(b)\}<\infty.$ Then there exists a probability $P(a,b)$ on $E$ with Laplace transform 
$L_{P(a,b)}(s)=L_{\mu}\{a(s+b)\}/L_{\mu}\{a(b)\}.$ This probability $P(a,b)$ is obtained in two steps: first take the image $\nu(dy)$ of $\mu(dx)$ 
by  the map $x\mapsto a^\top(x)=y$, where $a^\top$ is the adjoint of $a.$ Its Laplace transform is $L_{\nu}(s)=L_{\mu}\{a(s)\}.$ The second step constructs $P(a,b)(dy)$ as the probability $e^{-\<b,y\>}\nu(dy)/L_{\mu}\{a(b)\}$: it is a member of the exponential family generated by $\nu.$ 

Let us apply this to the case where $E$ is the Euclidean space of real symmetric matrices of order $d$ with scalar product $\<x,y\>=\tr(xy)$ and   where $\mu$ is $m(2p,k,d)$. Here  $D(\mu)=\mathcal{P}_d$. We take $b=(2\Sigma)^{-1}$ and  $a$ to be the linear transformation $s\mapsto a(s)= qsq^\top$, where $q$ is an invertible matrix of order $d$ 
such that 
\begin{equation}
\label{exq}
2(2\Sigma )^{-1}w(2\Sigma )^{-1}=q^{-1}I(k,d)(q^\top)^{-1}\;
.\end{equation}
We have $a^\top(x)=q^\top xq.$ The  distribution  $P(a,b)$ is the non-central  Wishart $\mathcal{NCW}(2p,w,\Sigma)$ distribution since 
\begin{equation}\label{MNCW}\frac{L_{\mu}\{a(s+b)\}}{L_{\mu}\{ a(b)\} }=\frac{1}{\det (I_d+2\Sigma	s)^{p}}e^{-\tr (2s(I_d+2\Sigma	s)^{-1}w)}.
\end{equation}
Using  \eqref{exq}, we show \eqref{MNCW} as follows:
\begin{multline*}
\tr[[(q^\top)^{-1}\{s+(2\Sigma)^{-1}\}^{-1}q^{-1}-(q^\top)^{-1}(2\Sigma)q^{-1}] I(k,d)]\\
=\tr[[ \{ s+(2\Sigma)^{-1}\}^{-1}-2\Sigma ] q^{-1}I(k,d)(q^\top)^{-1}]=-\tr (2s \{ I_d+2\Sigma s\}^{-1}w). \hspace{1cm}
\end{multline*}
The only thing left to prove is the existence of $q$ satisfying \eqref{exq}. To see this, since the matrix $2(2\Sigma )^{-1}w(2\Sigma )^{-1}$ of  $\overline{\mathcal{P}}_d$  has rank $k,$  we write
$2(2\Sigma )^{-1}w(2\Sigma )^{-1}=u\Delta u^\top$, where $\Delta=\mathrm{diag}(0,\ldots,0,\lambda_1^2,\ldots,\lambda_k^2)$ with $\lambda_i> 0$ and where $u$ is an orthogonal matrix 
of order $d.$ Taking $q= \mathrm{diag}(1,\ldots,1,\lambda_1^{-1},\ldots,\lambda_k^{-1})\, u^\top$ provides a solution of \eqref{exq}. 

The proof of the converse follows similar lines. 

Suppose now that $m(2p,k,d)$ exists and let $w\in \overline{\mathcal{P}_d}$ of rank $k$. Then, from the first statement of the proposition,  $\mathcal{NCW}(2p,\epsilon w,\Sigma)$ exists for all $\epsilon>0$.  The characteristic functions of $\mathcal{NCW}(2p,\epsilon w,\Sigma)$ and of  $\mathcal{NCW}(2p,0,\Sigma)$ are easily deduced from \eqref{NCW} and \eqref{NCWstat1}. Using these characteristic functions and Paul L\'evy's continuity theorem, we deduce that $\mathcal{NCW}(2p,\epsilon w,\Sigma)$ converges weakly to $\mathcal{NCW}(2p,0,\Sigma)$ when $\epsilon \to 0$. Therefore we can claim that $\mathcal{NCW}(2p,0,\Sigma)$ also exists. From Gindikin's theorem, $p$ must be in $\Lambda_d$ defined by \eqref{GYN}. \hfill $\square$

\bigskip \noindent \textbf{Example 1.} When $0\leq k\leq n\leq d$ we can use the above principle for constructing $NCW[n,2I(k,d),I_d]$ from $m(n,k,d).$ We take $q=I_d$ and $b=I_d/2.$
Since $a$ is the identity, we have
\begin{equation}
\label{SS}
m(n,k,d)(dx)=2^{dn/2}e^{2k}e^{\tr x/2}NCW [ n,2I(k,d),I_d] (dx).
\end{equation}

\bigskip
The next three propositions reformulate known facts in the language of the measures $m(2p,k,d).$ 

\bigskip \noindent %Formerly 2.2
\textbf{Proposition 2.} 
{\it Let $n$ and $k$ be integers such that $0\leq n,k\leq d.$ The measure $m(n,k,d)$  exists for $0\leq k\leq n\leq d$. Furthermore, the measure $m(d-1,d,d)$ exists.} 

\bigskip \noindent \textbf{Remark 1.} The proof of  the existence of $m(d-1,d,d)$ given below 
is easy, but its explicit  computation is done in Section~\ref{sec:4}.

\bigskip \noindent \textit{Proof.} Formula \eqref{SS} provides an explicit form of $m(n,k,d)$. For $2p>d-1$ the probability $\mathcal{NCW}(2p,I_d,I_d)$ exists as proved in Proposition~2 of \cite{LM2008}. This implies that 
$$\lim_{p\downarrow (d-1)/2} \mathcal{NCW}(2p,I_d,I_d)=\mathcal{NCW}(d-1,I_d,I_d)$$ exists by considering the characteristic function of  $ \mathcal{NCW}(2p,I_d,I_d)$ deduced from \eqref{NCW} and L\'evy's continuity theorem. From Proposition~1, we have the result. \hfill $\square$ 

\bigskip \noindent 
% Formerly 2.3
\textbf{Proposition 3.} 
{\it Suppose $d\geq 3.$ If $m(d-2,d-1,d)$ does not exist then, for any pair of integers $(n,k)$  such that $0\leq n<k<d$, $m(n,k,d)$ does not exist either. If $m(d-2,d,d)$ does not exist then, for any integer $n$ such that $0\leq n\leq d-2$,  $m(n,d,d)$ does not exist  either.}

\bigskip \noindent 
\textit{Proof.} Suppose that $m(n,k,d)$ exists for some pair $(n,k)$ such that $0\leq n<k<d.$  We define $m'(dx)$ as the measure on $\overline{\mathcal{P}_d}$ with Laplace transform  
$$
\int_{\overline{\mathcal{P}_d}}e^{-\tr(sx)}m'(dx)=\frac{1}{(\det s)^{(d-n-2)/{2}}} \, e^{\tr [s^{-1} \{I(d-1,d)-I(k,d)\}]}.
$$
Since the rank of $I(d-1,d) -I(k,d)$ is equal to $d-1-k$ and is less than or equal to $ d-n-2$ 
then, by Propositions~1--2,  $m'$ exists. Now we write the convolution
$m(n,k,d)*m'=m(d-2,d-1,d)$, which contradicts the non-existence of $m(d-2,d-1,d).$  Similarly, suppose that $m(d-2,d,d)$ does not exist and that there exists $n$ such that $0\leq n\leq d-2$ and such that $m(n,d,d)$ exists. Then $m(n,d,d)*m(d-2-n,0,d)=m(d-2,d,d)$ also leads to a contradiction. \hfill $\square$

\bigskip
An important result which is a consequence of     \cite{Mayerhofer}  is the following.

\bigskip \noindent 
% Formerly 2.4
\textbf{Proposition 4.} 
{\it If $d\geq 3$, the measure $m(d-2,d,d)$ does not exist.} 

\bigskip
 Apart from Proposition~9 below, our main result is as follows.

\bigskip \noindent 
% Formerly 2.5
\textbf{Proposition 5.} 
{\it If $d\geq 3$, the measure $m(d-2,d-1,d)$ does not exist.}

\bigskip
The proof will be given in Section~\ref{sec:6}. 
In the remainder of the paper we develop the tools that lead us to this proof. They will also enable us to give another proof of Proposition~4. Let us emphasize the fact that Propositions~1--5 lead to  a necessary and sufficient condition of existence of the distribution $\mathcal{NCW}(2p,w,\Sigma).$  It is worthwhile  to make the following synthesis  of  Propositions~1--5 about the existence of $\mathcal{NCW}(2p,w,\Sigma)$ in the following proposition.

\bigskip \noindent 
% Formerly 2.6
\textbf{Proposition 6.} 
{\it Let $\Sigma\in \mathcal{P}_d$ , $w\in \overline{\mathcal{P}_d}$ with rank $k \in \{ 0,\ldots,d\}$ and $p>0.$ Then the non-central Wishart distribution  $\mathcal{NCW}(2p,w,\Sigma)$ exists if and only either $2p\geq d-1$ or $2p=n\in \{0,1,\ldots,d-2\}$ and $0\leq k\leq n$. In particular, for $d=2$, the probability $\mathcal{NCW}(2p,w,\Sigma)$ exists if and only if $2p\geq 1$. For  an arbitrary $d$,  $\mathcal{NCW}(0,w,\Sigma)$ exists if  and only if $w=0$.} 

\bigskip \noindent 
\textit{Proof.} From Proposition~1, the existence of $\mathcal{NCW}(2p,w,\Sigma)$  is equivalent to the existence of $m(2p,k,d)$ when rank $w=k.$

\medskip
\noindent 
$\Rightarrow$: We have seen in Proposition~1 that $p$ must be in $ \Lambda_d$ as defined in \eqref{GYN}. If furthermore $2p = n\in \{0,\ldots,d-2\}$ let us show that  $0\leq k\leq n.$ Suppose the contrary $0\leq n<k$ holds. A reformulation of the first part of Proposition~3 is the following: if there exists $(n,k)$ such that $0\leq n<k<d$ and such that $m(n,k,d)$ exists then $m(d-2,d-1,d)$ exists. This contradicts the statement of Proposition~5. Thus the `if' part of Proposition~6 is proved. 

\medskip
\noindent
$\Leftarrow$: If $2p> d-1$ Proposition 2.2 of \cite{LM2008} proves the existence of $m(2p,k,d)$ without constraints on $k.$ Passing to the limit when $2p=d-1$ shows the existence of $m(d-1,k,d)$ also for any $k$. If $2p=n\leq d-2$ and if $0\leq k\leq n$, Proposition~2  shows that $m(n,k,d)$ exists. \hfill $\square$

\section{Computation of $m(1,2,2)$} 
\label{sec:3}

 In this section, we compute $m(1,2,2)$, which exists, as we know from Proposition~2. We will use only elementary tools.  We  parameterize the cone $\overline{\mathcal{P}_2}$ by the cone of revolution 
 $$
 C= \left\{(x,y,z)\in \mathbb{R}^3: x\geq \sqrt{y^2+z^2}\right\}
 $$ 
 using  the mapping $\varphi$ from $C$ to $\overline{\mathcal{P}_2}$ defined by
 \begin{equation}
 \label{VP}
 (x,y,z)\mapsto \varphi(x,y,z)= \left[\begin{array}{cc}x+y&z\\z&x-y\end{array}\right].\end{equation} Note  that  $\tr \{ \varphi (a,b,c)\varphi (x,y,z)\} =2ax+2by+2cz.$ 

%%% Formerly Proposition 3.1
\bigskip \noindent 
\textbf{Proposition 7.} 
{\it Consider the positive measure $\mu$ on $C$ such that, 
   for $a>\sqrt{b^2+c^2}$, we have 
$$
\frac{1}{\sqrt{a^2-b^2-c^2}} \, e^{\frac{2a}{a^2-b^2-c^2}}=\int_{C} \, e^{-2ax-2by-2cz}\mu(dx,dy,dz),
$$
i.e., such that the image of $\mu$ by $\varphi$ is $m(1,2,2).$ Then $\mu(dx,dy,dz)=r(dx,dy,dz)+f(x,y,z)\mathbf{1}_{C}(x,y,z)dxdydz$, where the singular part $r$  is the image of 
the measure $g(2\sqrt{y^2+z^2})dydz$ on $\mathbb{R}^2$ by the map $(y,z)\mapsto (x,y,z)=(\sqrt{y^2+z^2},y,z)$ with $g(t)=\{2 \cosh(2\sqrt{t})\}/(\pi t)$ and where, for $(x,y,z)\in C$,}
$$
f(x,y,z)=\frac{2}{\sqrt{\pi}}\sum_{k=0}^{\infty}\frac{(x^2-y^2-z^2)^{k}}{k!(k+1)!}\sum_{m=0}^{\infty} \frac{1}{\Gamma(m+2k+ {5}/{2})}\frac{(2x)^{m}}{m!}.
$$

\noindent 
\textit{Proof.} Using the notation  $D= {\partial}/{\partial x},$ the  Fa\`a di Bruno  differentiation formula states that, if $f(t)$ and $g(x)$ 
are functions with enough derivatives, then 
\begin{equation}\label{FdB}
D^nf \{ g(x)\} =\sum\frac{n!}{k_1!\cdots k_n!}(D^kf)\{ g(x)\}
\left\{ \frac{Dg(x)}{1!}\right\}^{k_1}\cdots
\left\{ \frac{D^ng(x)}{n!}\right\}^{k_n},
\end{equation}
where $k=k_1+\cdots+k_n$ and where the sum is taken on all integers 
 $k_j\geq 0$ such that  $k_1+2k_2+\cdots+nk_n=n.$ 
 For a reference, see, e.g., \cite{Roman}. We apply \eqref{FdB} to $g$ defined by $x\mapsto x^2-y^2-z^2$ for fixed $y,z$ and  to $f(t)=t^n.$ Noting that $D^3g=0$, we obtain
 \begin{equation}
 \label{FAA}
 \frac{\partial^n}{\partial x^n} \, (x^2-y^2-z^2)^n=n!^2\sum_{k_2=0}^{[n/2]}\frac{1}{k_2!}\times \frac{(x^2-y^2-z^2)^{k_2}}{k_2!}\times \frac{(2x)^{n-2k_2}}{(n-2k_2)!}.\end{equation} For simplification in the sequel we write
$E=e^{-2ax-2by-2cz}$ and $F=e^{-2a\sqrt{y^2+z^2}-2by-2cz}$. We now recall (see Formula~3.24 in \cite{LW}) that when $p>1/2$ we have for $a>\sqrt{b^2+c^2}$
 \begin{equation}\label{OUF2}\frac{1}{(a^2-b^2-c^2)^{p}}=\frac{2}{\sqrt{\pi}}\times\frac{1}{\Gamma(p) \Gamma(p- {1}/{2})}\int_{C}(x^2-y^2-z^2)^{p-{3}/{2}}Edx dy dz.
 \end{equation}
 Define 
\begin{eqnarray*}
I_k(n)=(2a)^{k}\int_{C}(x^2-y^2-z^2)^{n}Edx dy dz =\frac{\sqrt{\pi}}{2} \, n!\, \Gamma(n+ {3}/{2})\times\frac{(2a)^k}{(a^2-b^2-c^2)^{n+ {3}/{2}}} \, ,
\end{eqnarray*}
where we apply  \eqref{OUF2} for $p=n+ {3}/{2}$. The idea of the proof is to write the Laplace transform of $\mu$ as follows: 
\begin{multline}
\frac{e^{\frac{2a}{a^2-b^2-c^2}}}{\sqrt{a^2-b^2-c^2}} = \frac{1}{\sqrt{a^2-b^2-c^2}}+\sum_{n=0}^{\infty}\frac{(2a)^{n+1}}{(n+1)!} \frac{1}{(a^2-b^2-c^2)^{n+ {3}/{2}}} \\=\frac{1}{\sqrt{a^2-b^2-c^2}}+\frac{2}{\sqrt{\pi}}\sum_{n=0}^{\infty}\frac{1}{(n+1)!n!\Gamma(n+ {3}/{2}) }I_{n+1} (n)\;.\label{BF2}
\end{multline}
 A first step is to observe that 
 for  $k \in \{0,\ldots,n\}$, we have  
\begin{equation}
\label{OUF3} 
I_k(n)=\int_{C}\frac{\partial^k}{\partial x^k} \, (x^2-y^2-z^2)^nEdxdydz.\end{equation} 

Let us prove it by induction on $k$. It is true for $k=0$. Suppose that it is true for $k<n$ and let us show that \eqref{OUF3} is true for $k+1.$ 
Observe that for fixed $(y,z)$, the root $\sqrt{y^2+z^2}$ of the polynomial $x\mapsto (x^2-y^2-z^2)^n$ has order $n$ and this implies that $ {\partial^k}(x^2-y^2-z^2)^n/{\partial x^k}$ is zero for $x=\sqrt{y^2+z^2}.$ Using this remark and  integrating by parts with $V(x)=e^{-2ax}$ and $U(x)= {\partial^k}(x^2-y^2-z^2)^n / {\partial x^k}$, we compute the following integral:
\begin{equation}\label{INTP}
\int_{\sqrt{y^2+z^2}}^{\infty}2ae^{-2ax}\frac{\partial^k}{\partial x^k}(x^2-y^2-z^2)^ndx=\int_{\sqrt{y^2+z^2}}^{\infty}e^{-2ax}\frac{\partial^{k+1}}{\partial x^{k+1}} \, (x^2-y^2-z^2)^ndx.
\end{equation} 
With \eqref{INTP} we are in position to prove  \eqref{OUF3}. We have
\begin{eqnarray*}
I_{k+1}(n)&=&2a\int_{C}\frac{\partial^k}{\partial x^k} \, (x^2-y^2-z^2)^nEdxdydz =
\int_{\mathbb{R}^2}e^{-2by-2cz}\left[\int_{\sqrt{y^2+z^2}}^{\infty}2ae^{-2ax}\frac{\partial^k}{\partial x^k} \,
(x^2-y^2-z^2)^ndx\right]dydz
\\&=&\int_{\mathbb{R}^2}e^{-2by-2cz}\left[\int_{\sqrt{y^2+z^2}}^{\infty}e^{-2ax}
\frac{\partial^{k+1}}{\partial x^{k+1}} \, (x^2-y^2-z^2)^ndx\right]dydz
=\int_{C}\frac{\partial^{k+1}}{\partial x^{k+1}} \, (x^2-y^2-z^2)^nEdxdydz,
\end{eqnarray*}
which proves \eqref{OUF3}. We will need \eqref{OUF3} only for $k=n$. 

The second step is to express $I_{n+1}(n)$  as the Laplace transform of a positive measure.  We compute $I_n(n)$ as expressed \eqref{OUF3} by using  again an integration by parts. The new fact for $k=n$ is that the integrated part will not disappear and will provide a term for  the singular measure $s$ given  in the statement of the theorem. This calculation of the integrated part will use \eqref{FAA}. Taking
$V(x)=-e^{-2ax}$ and $U(x)= {\partial^n}(x^2-y^2-z^2)^n / {\partial x^n}$, we write 
\begin{eqnarray}    
I_{n+1}(n)=2a  I_n(n)&=&\int_{\mathbb{R}^2}e^{-2by-2cz}\left[\int_{\sqrt{y^2+z^2}}^{\infty}2ae^{-2ax}\frac{\partial^n}{\partial x^n} \, (x^2-y^2-z^2)^ndx\right]dydz =A_n+S_n
\end{eqnarray}
with
\begin{eqnarray}
A_n&=&\int_{C}\frac{\partial^{n+1}}{\partial x^{n+1}}(x^2-y^2-z^2)^nEdxdydz, \label{ABS}\\S_n&=&n!\int_{\mathbb{R}^2}e^{-2by-2cz}\left[-e^{-2ax}(2x)^n\right]_{\sqrt{y^2+z^2}}^{\infty}dydz = n!\int_{\mathbb{R}^2} \Bigl(2\sqrt{y^2+z^2} \Bigr) ^nFdydz,
\label{SING}
\end{eqnarray}
where  \eqref{SING} comes from \eqref{FAA} by keeping only the term $k_2=0.$ 
We will carry this value of  $I_{n+1}(n)=A_n+S_n$ in   \eqref{BF2}. Doing this, we can guess that $S_n$ will contribute to the singular part of $m(1,2,2).$ But the term $1/{\sqrt{a^2-b^2-c^2}}$ in \eqref{BF2} will also contribute to it. 

More specifically, the  third step  of the proof is to represent the function on $C\setminus \partial C$ defined by $(a,b,c)\mapsto  {1}/{\sqrt{a^2-b^2-c^2}}$ as a Laplace transform. Exploiting the Gaussian integral below, we obtain 
$$
\frac{1}{\sqrt{a^2-b^2-c^2}} =\frac{2}{\pi}\int_{\mathbb{R}^2}e^{-2a(u^2+v^2)-2b(u^2-v^2)-4cuv}dudv
= \frac{2}{\pi}\int_{\mathbb{R}^2} \Bigl(2\sqrt{y^2+z^2} \Bigr)^{-1}Fdydz.
$$
To obtain the second identity, observe that the map on $\{(u,v): u>0\}$ defined by $y=u^2-v^2,\ z=2uv$ is a bijection with $\mathbb{R}^2$; the same is true with  $\{(u,v): u<0\}$. Furthermore $dydz=4(u^2+v^2)dudv=4\sqrt{y^2+z^2}dudv$  and therefore $dudv = {dydz}/({4\sqrt{y^2+z^2})}$. This leads to the conclusion.

Now comes the fourth and final step. We use $\Gamma(n+ {1}/{2})= {(2n)!}\sqrt{\pi}/({4^n n!})$ and we consider the function 
$$
g(t)=\frac{2}{\pi t}+\frac{2}{\sqrt{\pi}}\sum_{n=0}^{\infty}\frac{1}{(n+1)!\Gamma(n+{3}/{2}) } \, t^n=\frac{2}{\pi t}\cosh(2\sqrt{t}).
$$
We then define the measure $r(dx,dy,dz)$ concentrated on the boundary 
$$
\partial C= \left\{(x,y,z): x=\sqrt{y^2+z^2} \right\}
$$ 
of the cone $C$ to be the image of 
the measure on $\mathbb{R}^2$ 
$$
g\Bigl( 2\sqrt{y^2+z^2} \Bigr)dydz=\frac{1}{\pi\sqrt{y^2+z^2}}\cosh\{2^{3/2}(y^2+z^2)^{1/4}\} dydz
$$
 by the map $(y,z)\mapsto (x,y,z)=(\sqrt{y^2+z^2},y,z).$ This measure $r$ will be the singular part of the image $\mu$ of $m(1,2,2)$ by the reciprocal of $\varphi$ defined by \eqref{VP}:  
\begin{eqnarray}
\nonumber\int_C Eds&=&\int_{\R^2}      g \Bigl (2\sqrt{y^2+z^2} \Bigr)Fdydz =\nonumber\int_{\R^2}F\frac{dydz}{\pi\sqrt{y^2+z^2}}+\frac{2}{\sqrt{\pi}}\sum_{n=0}^{\infty}\int_{\R^2}\frac{2^n(\sqrt{y^2+z^2})^n}{(n+1)!\Gamma(n+ {3}/{2})}Fdydz\\ \label{SING2}
&=&\frac{1}{\sqrt{a^2-b^2-c^2}}+\frac{2}{\sqrt{\pi}}\sum_{n=0}^{\infty}\frac{S_n}{(n+1)!n!\Gamma(n+ {3}/{2})} \, .
\end{eqnarray}

Finally we focus on the absolutely continuous part of $\mu$. We will need the following formula, similar to \eqref{FAA} and also obtained using  \eqref{FdB}: 
$$
 \frac{\partial^n}{\partial x^n}(x^2-y^2-z^2)^{n-1}=n!\, (n-1)! \, \sum_{k_2=1}^{[n/2]}\frac{1}{(k_2-1)!}\times \frac{(x^2-y^2-z^2)^{k_2-1}}{k_2!}\times \frac{(2x)^{n-2k_2}}{(n-2k_2)!}\, .
 $$
The absolutely continuous part of  $\mu$  is given by \eqref{ABS}  and \eqref{BF2}. Its  density is 
\begin{eqnarray}
f(x,y,z)&=&\label{ABS2}
\frac{2}{\sqrt{\pi}}\sum_{n=0}^{\infty}\frac{1}{(n+1)! \, n! \, \Gamma(n+ {3}/{2}) }\frac{\partial^{n+1}}{\partial x^{n+1}} \, (x^2-y^2-z^2)^n\\
&=&\nonumber \frac{2}{\sqrt{\pi}}\sum_{n=2}^{\infty} \frac{1}{(n-1)! \, n! \, \Gamma(n+ {1}/{2}) }\frac{\partial^{n}}{\partial x^{n}} \, (x^2-y^2-z^2)^{n-1}
\\&=&\nonumber \frac{2}{\sqrt{\pi}}\sum_{n=2}^{\infty}\frac{1}{\Gamma(n+ {1}/{2}) }
\sum_{k_2=1}^{[n/2]}\frac{1}{(k_2-1)!}\times \frac{(x^2-y^2-z^2)^{k_2-1}}{k_2!}\times \frac{(2x)^{n-2k_2}}{(n-2k_2)!}
\\&=& 
\frac{2}{\sqrt{\pi}}\sum_{k=0}^{\infty}\frac{(x^2-y^2-z^2)^{k}}{k! \, (k+1)!}\sum_{m=0}^{\infty} \frac{1}{\Gamma(m+2k+ {5}/{2})}\frac{(2x)^{m}}{m!}\,. \notag
\end{eqnarray}
From \eqref{ABS2}, \eqref{BF2} and \eqref{ABS} the Laplace transform of $f$ is 
\begin{equation}\label{ABS4}
\int_C E fdxdydz=\frac{2}{\sqrt{\pi}}\sum_{n=0}^{\infty}\frac{A_n}{(n+1)! \, n! \, \Gamma(n+{3}/{2})} \, .
\end{equation}
Let us add \eqref{SING2} and \eqref{ABS4} and use \eqref{BF2}. What we get shows that the parameterization $\mu$ of $m(1,2,2)$ by $\varphi$ is the sum of $r$ and of the absolutely continuous part with density $f$. Formula {\eqref{ABS4}} shows that $f$ is as given in Proposition~7. \hfill
$\square$

\bigskip \noindent 
\textbf{Remark 2.} The following point  is essential for  understanding  Section 4. Remark first that the image by $\varphi$ of the measure $r(dx,dy,dz)$ is concentrated on the set $S_1\subset \overline{\mathcal{P}}_2$ of matrices of rank $1$. Any element of $S_1$ can be written as 
$$
u\left[\begin{array}{cc}\lambda_1&0\\0&0\end{array}\right] u^\top,
$$ 
where $u$ is an orthogonal matrix of $\mathbb{O}(2)$ and $\lambda_1>0.$ We can compute the image of $r(dx,dy,dz)$  by the map 
\begin{equation}\label{MT}
\left[\begin{array}{cc}x+y&z\\z&x-y\end{array}\right]=\left[\begin{array}{cc}\sqrt{y^2+z^2}+y&z\\z&\sqrt{y^2+z^2}-y\end{array}\right]\mapsto \lambda_1=2\sqrt{y^2+z^2}.\end{equation} 
If $A_t=\{(x,y,z) :   2\sqrt{y^2+z^2}<t\}$, then using polar coordinates $y=\lambda_1\cos\alpha$, $z=\lambda_1\sin \alpha$ with Jacobian equal to $ {\lambda_1}/{2}$, we have
$$r(A_t)=\int_{A_t}r(dx,dy,dz)=\int_{2\sqrt{y^2+z^2}<t}g \Bigl (2\sqrt{y^2+z^2} \Bigr)dydz=\frac{\pi}{2}\int_0^tg(\lambda_1)\lambda_1d\lambda_1.$$  
Since $g(\lambda_1)= {2}(\pi \lambda_1)^{-1}\cosh 2\sqrt{\lambda_1}$ , the image of the measure $r$ by the map \eqref{MT} is 
\begin{equation}\label{M1}
\cosh (2\sqrt{\lambda_1})\textbf{1}_{(0,\infty)}(\lambda_1)d\lambda_1.
\end{equation} 
Now an important observation is the following: consider the measure $m(1,1,1)(d\lambda)$ on $(0,\infty)$ whose Laplace transform is, for $s>0,$
$$
\frac{1}{\sqrt{s}} \, e^{1/s}=\sum_{n=0}^{\infty}\frac{1}{n!s^{n+ {1}/{2}}}=\int_0^{\infty}e^{-s\lambda}\sum_{n=0}^{\infty}\frac{\lambda^{n-{1}/{2}}}{n!\Gamma(n+{1}/{2})}d\lambda=\frac{1}{\sqrt{\pi}}\int_0^{\infty}e^{-s\lambda}\frac{1}{\sqrt{\lambda}}\cosh(2\sqrt{\lambda})d\lambda.
$$ 
 This last line implies 
\begin{equation}
\label{m111}
m(1,1,1)(d\lambda)=\frac{1}{\sqrt{\pi}}\frac{1}{\sqrt{\lambda}}\cosh(2\sqrt{\lambda})\textbf{1}_{(0,\infty)}(\lambda)d\lambda
\end{equation}
 and one observes that $m(1,1,1)$ is quite close to \eqref{M1}. To summarize this remark, the singular part of $m(1,2,2)$ can be seen as the image of $\sqrt{\pi \lambda_1}m(1,1,1)(d\lambda_1)\otimes du$ by the map 
 $$
 (u,\lambda_1)\mapsto u\left[\begin{array}{cc}\lambda_1&0\\0&0\end{array}\right] u^\top
 $$ 
 from $(0,\infty)\times \mathbb{O}(2) $, where $du$ is the uniform probability on $\mathbb{O}(2).$ This is the key to the generalization  of the computation from  $m(1,2,2)$ to the computation of  $m(d-1,d,d)$ for $d\geq 2$ done in the next section.  

\section{Computation of the measure $m(d-1,d,d)$} 
\label{sec:4}

Before stating Proposition~9 which describes $m(d-1,d,d)$ we set some notations,   we recall some facts about zonal  functions and polynomials and we prove three lemmas. The Lebesgue measure $dx$ on the space of real symmetric matrices of order $d$ has the normalization  associated to the Euclidean structure  given by $\<x,y\>=\tr(xy)$. Note that  Muirhead \cite{Muirhead}  has a different normalization. As mentioned on page ix of the Introduction of \cite{Muirhead}, zonal functions are the essential tool of non-central  distribution theory. Our reference will be Faraut and Koranyi \cite {FK}, abbreviated from now on as FK. 

\subsection{Zonal functions} 
\label{sec:4.1}

Let $\mathcal{E}_d$ be the set of sequences  $\kappa=(m_1,\ldots,m_d)$ of $d$ integers such that $m_1\geq  \cdots\geq m_d\geq 0$. For $\kappa\in \mathcal{E}_d$, we consider the two zonal polynomials 
$$
\Phi^{(d)}_{\kappa}(x)=\Phi^{(d)}_{m_1,\ldots,m_d}(x),\ \ C^{(d)}_{\kappa}(x)=C^{(d)}_{\kappa}(I_d)\Phi^{(d)}_{\kappa}(x),
$$ 
where $C^{(d)}_{\kappa}(I_d)$ is defined below in \eqref{CST}. In FK p.~228, the $\Phi_{\kappa}$ are  called spherical  rather than zonal polynomials, and on p.~234 the notation $Z_{\kappa}$ is used instead of our notation  $C^{(d)}_{\kappa}$. We use the definitions given in FK, while \cite{Muirhead}  and \cite{Takemura}  have other ways to introduce the zonal polynomials. To define $\Phi^{(d)}_{\kappa}$ we consider 
$$
\Delta_{\kappa}(x)=\Delta_1(x)^{m_1-m_2}\Delta_2(x)^{m_2-m_3}\cdots\Delta_{d-1}(x)^{m_{d-1}-m_d}\Delta_d(x)^{m_d},
$$
where for  the real symmetric matrix  $x=(x_{ij})_{1\leq i,j\leq d}$, the function  $\Delta_k(x)=\det(x_{ij})_{1\leq i,j\leq k}$ is the principal determinant of $x$ of order $k.$ The function $\Phi^{(d)}_{\kappa}$ is then defined  by 
\begin{equation}\label{Phi}
\Phi^{(d)}_{\kappa}(x)=\int_{\mathbb{O}(d)}\Delta_{\kappa}(uxu^\top)du,\end{equation}
where $du$ is the Haar probability on the orthogonal group $\mathbb{O}(d).$ 
When $x\in \mathcal{P}_d$  definition  \eqref{Phi} makes sense even when 
$m_1,\ldots,m_d$ are complex numbers. In that case $\Phi^{(d)}_{m_1,\ldots,m_d}(x)$ is no longer a polynomial and is called a zonal function. 

From \eqref{Phi} one obtains easily the following formula for any $x$ and $y$ in $\mathcal{S}_d$:
$$
\int_{\mathbb{O}(d)} \Phi^{(d)}_{\kappa}(xuyu^\top)du=\Phi^{(d)}_{\kappa}(x)\Phi^{(d)}_{\kappa}(y).
$$
To give the value of  the  constant $C^{(d)}_{\kappa}(I_d)$, we need the notations $\ell(\kappa)=\max\{j: m_j>0\}$, $|\kappa|=m_1+\cdots+m_d$ and  
$$
\Gamma_d( z_1,\ldots,z_d)=\prod_{j=1}^d\Gamma \left(z_j-\frac{j-1}{2} \right)
$$
defined whenever $z_j- ({j-1})/{2}>0$ for all $j \in \{ 1,\ldots,d\}$. If $p$ is a real number, we use the notational convention 
$$
\Gamma_d( z+p) = \Gamma_d( z_1+p,\ldots,z_d+p).
$$

If $\kappa\in \mathcal{E}_d$ and $p>(d-1)/2$, we define the Pochhammer symbol as $
(p)_{\kappa}=\Gamma_d( \kappa+p)/\Gamma_d( p)$.
Since $p\mapsto (p)_{\kappa}$ is a polynomial, we extend its definition to the whole line. 
If $\kappa\in \mathcal{E}_d$, the constant $C^{(d)}_{\kappa}(I_d)$ is 
\begin{eqnarray}\label{CST}
C^{(d)}_{\kappa}(I_d)&=&C^{(d)}_{m_1,\ldots,m_d}(I_d)=2^{2|\kappa|}|\kappa|!\left(\frac{d}{2}\right)_{\kappa}\frac{\prod_{1\leq i<j\leq \ell(\kappa)}(2m_i-2m_j-i+j)}{\prod_{i=1}^{\ell(\kappa)}\{2m_i+\ell(\kappa)-i\}!}\\&=&\frac{|\kappa|!}{((d+1)/2)_{\kappa}}\times\prod_{1\leq i\leq j\leq d}\frac{B[(j-i+1)/2, {1}/{2}]}{B[m_i-m_j+(j-i+1)/2,{1}/{2}]}.\label{CST2}
\end{eqnarray}

Form \eqref{CST} of the spherical polynomials is given in \cite{ Muirhead} p.~237 formula (38), where there is a reference to \cite{Constantine} for a proof. Form  \eqref{CST2}  can be  proved from FK by combining Propositions XI. 4.2 (i), p.~230, and XI.4.4, p.~232, and the definition of  $Z_{\kappa}(x)=C_{\kappa}^{(d)}(x)$ on the last line of p.~234. We never consider $C^{(d)}_{\kappa}(x)$ if $\kappa\notin \mathcal{E}_d.$ The exact value of $C^{(d)}_{\kappa}(I_d)$ will be crucial  in the proof of Proposition~9 when we  need  the following formula  which is Formula (3),  p.~259, of Muirhead \cite{Muirhead}: 
 \begin{equation}\label{EXPTR}e^{\tr x}=\sum_{\kappa\in \mathcal{E}_{d}}\frac{1}{|\kappa|!}C^{(d}_{\kappa}(x)\;.\end{equation} 
Since we have introduced the zonal polynomials, we use them here to recall   the result of  \cite{Peddada} 
 leading to  the  elegant  proof by Peddada and Richards of Gindikin's Theorem.

%%% Former Proposition 4.1
\bigskip \noindent 
\textbf{Proposition~8.} 
{\it Let $p>0$ and let $W$ be a random variable of $\overline{\mathcal{P}_d}$ such that $\esp(e^{\tr sW})= \{\det(I_d-s)\}^{-p}$ when $I_d-s$ is positive definite. Then, for all $\kappa\in \mathcal{E}_d$, we have $\esp \{ \Phi_{\kappa}^{(d)}(W)\}=(p)_{\kappa}$ and $W$ exists if and only if $p$ is in $\Lambda_d$ as defined by~\eqref{GYN}.}

\bigskip \noindent \textit{Proof.} Consider the scalar product  for polynomials on $\mathcal{S}_d$ defined  in FK, p.~220, as
$$
\langle p,q\rangle=p\left(\frac{\partial}{\partial x}\right)q(x)\Big|_{x=0}.
$$
Then (see FK, p.~234) the spherical polynomials are an orthogonal family. In particular 
\begin{eqnarray}\nonumber
\esp \{\Phi^{(d)}_{\kappa}(W) \}&=&\Phi^{(d)}_{\kappa}\left(\frac{\partial}{\partial s}\right)\esp (e^{\tr sW} )\big|_{s=0}=\Phi^{(d)}_{\kappa}\left(\frac{\partial}{\partial s}\right)\{\det(I_d-s)\}^{-p}\big|_{s=0}\\\label{PRFK3}&=&\Phi^{(d)}_{\kappa}\left(\frac{\partial}{\partial s}\right)\sum_{\kappa'\in \mathcal{E}_d}\frac{(p)_{\kappa}}{|\kappa'|!}C^{(d)}_{\kappa'}(s)\big|_{s=0}\\\label{PRFK4}&=&(p)_{\kappa} \|\Phi^{(d)}_{\kappa}\|^2C^{(d)}_{\kappa}(I_d)=(p)_{\kappa} \end{eqnarray}
where \eqref{PRFK3} comes from \eqref{EXPTR}, \eqref{PRFK4} comes from the orthogonality of the spherical functions, and  the norm of $\Phi^{(d)}_{\kappa}$ is given in FK, p.~234. Since $\Phi_{\kappa}^{(d)}(W)\geq 0$ we have $(p)_{\kappa}\geq 0$ for all $\kappa\in \mathcal{E}_d$ and this implies that $p\in \Lambda_d$. \hfill $\square$ 

\bigskip
It is interesting to compute the zonal polynomials for $d=2.$  This is given as Exercise~5, p.~237, in FK. 
Define the Legendre polynomials $(P_k)_{k=0}^{\infty}$ by their generating function 
$$
\sum_{k=0}^{\infty}P_k(x)z^k=\frac{1}{\sqrt{1-2zx+z^2}}.
$$
Further let 
$$
x=\left[\begin{array}{cc}a+b&c\\c&a-b\end{array}\right]
$$ 
in  $\mathcal{P}_2.$ Then for $(m_1,m_2)\in \mathcal{E}_2$ we have 
\begin{equation}\label{LP2}
\Phi_{m_1,m_2}(x)=(a^2-b^2-c^2)^{(m_1+m_2)/2}P_{m_1-m_2}
\left(\frac{a}{\sqrt{a^2-b^2-c^2}} \right).
\end{equation}
We now detail the proof of \eqref{LP2}.  The Legendre polynomial $P_k$ satisfies 
\begin{equation}\label{LPI}
P_k(\cosh t)=\frac{1}{\pi}\int_0^{\pi}(\cosh t+\cos u\sinh t)^kdu.
\end{equation}

To check this, call $Q_k$ the right-hand side of   \eqref{LPI}. The computation of  $\sum_{k=0}^{\infty}Q_kz^k$ gives ${(1-2zx+z^2)^{-1/2}}$ for $|z|$ small enough. This proves \eqref{LPI}.
Let 
$$ 
R(\theta)=\left[\begin{array}{cc}\cos \theta&-\sin \theta\\\sin \theta&\cos \theta\end{array}\right],\quad 
J=\left[\begin{array}{cc}-1&0\\0&1\end{array}\right]
$$ 
and observe that
$\mathbb{SO}(2)=\{R(\theta): \theta\in \R\}$ and $ \mathbb{O}(2)\setminus \mathbb{SO}(2)=J\mathbb{SO}(2)$.
Writing
$$
\left(\begin{array}{c}B\\C\end{array}\right)=R(2\theta)\left(\begin{array}{c}b\\c\end{array}\right),  
$$
a small  calculation yields
\begin{eqnarray}\label{RT1}R(\theta)\left[\begin{array}{cc}a+b&c\\c&a-b\end{array}\right]R(-\theta)&=&\left[\begin{array}{cc}a+B&C\\C&a-B\end{array}\right], \\\label{RT2} JR(\theta)\left[\begin{array}{cc}a+b&c\\c&a-b\end{array}\right]R(-\theta)J&=&\left[\begin{array}{cc}a+B&-C\\-C&a-B\end{array}\right]\;.\end{eqnarray}
The two formulas \eqref{RT1} and \eqref{RT2} enable us to compute the zonal polynomial
\begin{eqnarray*}
\Phi_{m_1,m_2}(x)&=&\int_{\mathbb{O}(2)}\Delta_{m_1,m_2}(uxu^\top)du =
\frac{1}{2}\int_{\mathbb{SO}(2)}\Delta_{m_1,m_2}(uxu^\top)du+\frac{1}{2}\int_{\mathbb{O}(2)\setminus\mathbb{SO}(2)}\Delta_{m_1,m_2}(uxu^\top)du\\
&=& \frac{1}{4\pi}\int_0^{2\pi} [ \Delta_{m_1,m_2}\{ R(\theta)xR(-\theta)\} +\Delta_{m_1,m_2} \{ JR(\theta)xR(-\theta)J \} ] d\theta\\
&=&(a^2-b^2-c^2)^{m_2}\frac{1}{2\pi}\int_0^{2\pi}(a+B)^{m_1-m_2}d\theta\\&=&(a^2-b^2-c^2)^{m_2}\frac{1}{\pi}\int_0^{\pi}(a+\sqrt{b^2+c^2}\cos \theta)^{m_1-m_2}d\theta\\&=&(a^2-b^2-c^2)^{\frac{m_1+m_2}{2}}\frac{1}{\pi}\int_0^{\pi}\left(\frac{a}{\sqrt{a^2-b^2-c^2}}+\frac{\sqrt{b^2+c^2}}
{\sqrt{a^2-b^2-c^2}}\cos \theta\right)^{m_1-m_2}d\theta\;.
\end{eqnarray*} Using form  \eqref{LPI} of the Legendre polynomial yields  \eqref{LP2}. Using \eqref{CST2} we also obtain
\begin{equation}\label{CON2}
C^{(2)}_{m_1,m_2}(I_2)=\frac{(m_1+m_2)!}{(m_1-m_2)! \, m_2!}\times\frac{1}{({3}/{2}+m_1-m_2)_{m_2}}.
\end{equation}

\subsection{Three properties of zonal functions}
\label{sec:4.2}

%%% Former Lemma 4.2
\bigskip \noindent \textbf{Lemma 1.} 
{\it Let $x=\left[\begin{array}{cc}x_1&x_{12}\\x_{21}&x_2\end{array}\right]\in \mathcal{P}_d$ and  $[x]_1=x_1\in \mathcal{P}_{d-1}$. Then for all complex numbers $m_1,\ldots,m_d$ we have} 
$$
\Phi^{(d)}_{m_1,\ldots,m_d}(x)=(\det x)^{m_d}\int_{\mathbb{O}(d)}\Phi^{(d-1)}_{m_1,\ldots,m_{d-1}}([uxu^\top]_1)du.
$$

\bigskip \noindent \textit{Proof.} Consider $v=\left[\begin{array}{cc}v_1&0\\0&1\end{array}\right]\in \mathbb{O}(d)$, where $v_1\in \mathbb{O}(d-1).$ Observe that for any $y\in  \mathcal{P}_d$ we have
\begin{equation}\label{Ja}
[vyv^\top]_1=v_1[y]_1v_1^\top.\end{equation}
We write 
\begin{eqnarray}\Phi^{(d)}_{m_1,\ldots,m_d}(x)&=&\int_{\mathbb{O}(d)}\Delta_{m_1,\ldots,m_{d-1},m_d}(uxu^\top)du \label{J0}
\\&=&(\det x)^{m_d}\int_{\mathbb{O}(d)}\Delta_{m_1,\ldots,m_{d-1}}([uxu^\top]_1)du \label{J1}
\\&=&(\det x)^{m_d}\int_{\mathbb{O}(d)}\Delta_{m_1,\ldots,m_{d-1}}([vuxu^\top v^\top]_1)du \label{J2}
\\&=&(\det x)^{m_d}\int_{\mathbb{O}(d)}\Delta_{m_1,\ldots,m_{d-1}}(v_1[uxu^\top v^\top]_1v_1^\top)du \label{J3}
\\&=&(\det x)^{m_d}\int_{\mathbb{O}(d)}\left[ \int_{\mathbb{O}(d-1)}\Delta_{m_1,\ldots,m_{d-1}} \{ v_1[uxu^\top v^\top]_1v_1^\top\} dv_1\right] du\label{J4}
\\&=&(\det x)^{m_d}\int_{\mathbb{O}(d)}\Phi^{(d-1)}_{m_1,\ldots,m_{d-1}}([uxu^\top]_1)du . \label{J5}\end{eqnarray}
In the equations above, \eqref{J0} comes from the definition of $\Phi^{(d)}_{\kappa}(x)$, \eqref{J1} separates the roles of $[uxu^\top]_1$ and $\det (uxu^\top)=\det x$
in the definition of $\Delta_{\kappa}(uxu^\top),$ \eqref{J2} uses the fact that $du$ is the Haar probability, \eqref{J3} follows from \eqref{Ja} applied to $y=uxu^\top$, \eqref{J4} uses the fact that the Haar measure $dv_1$ of $\mathbb{O}(d-1)$ has mass 1 and \eqref{J5} comes from the definition of $\Phi^{(d-1)}_{m_1,\ldots,m_{d-1}}(x)$. \hfill $\square$ 

%%% Former Lemma 4.3
\bigskip \noindent 
\textbf{Lemma 2.} 
{\it If $x \in \mathcal{P}_d$, then $\Phi^{(d)}_{m_1,\ldots,m_d}(x^{-1})=\Phi^{(d)}_{-m_d,\ldots,-m_1}(x)$.}

\bigskip \noindent \textit{Proof.}  Define $p\in \mathbb{O}(d)$ by 
$$
p=\left[\begin{array}{ccccc}
0&0&\ldots&0&1\\
0&0&\ldots&1&0\\
\vdots&\vdots&\vdots&\vdots&\vdots\\
1&0&\ldots&0&0\end{array}\right]
$$
and define $\Delta^\top_{m_1,\ldots,m_{d-1},m_d}(x)=\Delta_{m_1,\ldots,m_{d-1},m_d}(pxp^\top).$
We can write $\Phi^{(d)}_{\kappa}(x^{-1})$ as
\begin{eqnarray}\int_{\mathbb{O}(d)}\Delta_{m_1,\ldots,m_{d-1},m_d}(ux^{-1}u^\top)du 
&=&\int_{\mathbb{O}(d)}\Delta^\top_{-m_d,\ldots,-m_{1}}(uxu^\top)du \label{K1}
\\&=&\int_{\mathbb{O}(d)}\Delta_{-m_d,\ldots,-m_{1}}(uxu^\top)du \label{K2}
\end{eqnarray}
In this list \eqref{K1} comes  FK, Proposition VII.1.5 (ii), p.~127, which states that
$\Delta_{m_1,\ldots,m_{d-1},m_d}(x^{-1})= \Delta^\top_{-m_d,\ldots,-m_1}(x)$, 
and \eqref{K2} comes from the invariance of the Haar probability $du$ on $\mathbb{O}(d)$ by $u\mapsto pu.$\hfill $\square$

%%% Former Lemma 4.4
\bigskip \noindent \textbf{Lemma~3.} 
{\it If $x \in \mathcal{P}_d$, then $\Phi^{(d)}_{m_1,\ldots,m_d}(x)(\det x)^p=\Phi^{(d)}_{m_1+p,\ldots,m_d+p}(x)$.}

\bigskip \noindent \textit{Proof.} By definition, 
$$
\Phi^{(d)}_{\kappa}(x)(\det x)^p=\int_{\mathbb{O}(d)}\Delta_{m_1,\ldots,m_{d}}(uxu^\top)(\det uxu^\top)^pdu
=\int_{\mathbb{O}(d)}\Delta_{m_1+p,\ldots,m_{d}+p}(uxu^\top)du=\Phi^{(d)}_{m_1+p,\ldots,m_{d}+p}(x),
$$
which completes the argument. \hfill $\Box$

\subsection{The derivation of  $m(d-1,d,d)$}
\label{sec:4.3}

%%% Former Proposition 4.5
\bigskip \noindent \textbf{Proposition~9.} 
{\it Define the singular measure $r(dt)$ on $\overline{\mathcal{P}_d}$ and concentrated on the set $S_{d-1}$ of symmetric matrices of rank $d-1$ as  the image of  the product measure 
$$
\frac{(\pi \det x)^{1/2}}{\Gamma(d/2)}\, m(d-1,d-1,d-1)(dx)\otimes du 
$$ 
by the map
from $\mathcal{P}_{d-1}\times \mathbb{O}(d)$ to $\overline{\mathcal{P}_{d}}$ defined by  
$$
(x,u)\mapsto t=u\left[\begin{array}{cc}x&0\\0&0\end{array}\right]u^\top=u\tilde{x}u^\top.
$$
 Define 
 \begin{equation}
 \label{ABS6}
 f_d(t)=(\det t)^{-1}\left[ \sum_{\kappa\in \mathcal{E'}_d}\frac{C^{(d)}_{\kappa}(t)}{|\kappa|! \, \Gamma_d \{ \kappa+ {(d-1)}/{2} \}}\right] ,
 \end{equation}
where $\mathcal{E'}_d=\{\kappa\in \mathcal{E}_d\ ;\ m_d>0\}.$ Then 
$m(d-1,d,d)(dt)=r(dt)+f_d(t)\textbf{1}_{\mathcal{P}_d}(t)dt$.} 

\bigskip \noindent  \textit{Proof.} The function $f_d(t)$ is an analytic function around $t=0$ because, from the definition \eqref{Phi} of $\Phi^{(d)}_{\kappa}$, the polynomial $C^{(d)}_{m_1,\ldots,m_d}(t)$ is divisible by $(\det t)^{m_d}.$ Therefore $(\det t)^{-1}C^{(d)}_{\kappa}(t)$ is a polynomial when $\kappa\in \mathcal{E}'_d.$ 
Recall the following (see FK, Lemma XI.2.3, p.~226, or Muirhead \cite{Muirhead} Theorem 7.2.7, p.~248):
\begin{equation}
\label{BZ}
\int_{\mathcal{P}_{d}}e^{-\tr(sx)}(\det x)^{p-(d+1)/2}\frac{\Phi^{(d)}_{\kappa}(x)}{\Gamma_d(\kappa+p)} \, dx=\Phi^{(d)}_{\kappa}(s^{-1})(\det s)^{-p}.
\end{equation}
Note that the choice of the suitable  Lebesgue measure $dx$ is crucial in \eqref{BZ}.
Formula \eqref{BZ} holds true for $p+m_d>(d-1)/2.$ For $m_d=0$, this was proven in the references  \cite{Muirhead} and FK. When $m_d>0$ we observe that 
\begin{equation}\label{PS1}\Phi^{(d)}_{\kappa}(x)=\Phi^{(d)}_{m_1,\ldots,m_d}(x)=(\det x)^{m_d}\Phi^{(d)}_{m_1-m_d,\ldots,m_{d-1}-m_d,0}(x)=(\det x)^{m_d}\Phi^{(d)}_{\kappa-m_d}(x).\end{equation}
As a consequence 
\begin{eqnarray}
\int_{\mathcal{P}_{d}}e^{-\tr(sx)}(\det x)^{p-(d+1)/2}\frac{\Phi^{(d)}_{\kappa}(x)}{\Gamma_d(\kappa+p)}dx  \label{PS2}&=&\int_{\mathcal{P}_{d}}e^{-\tr(sx)}(\det x)^{p+m_d-(d+1)/2}\frac{\Phi^{(d)}_{\kappa-m_d}(x)}{\Gamma_d(\kappa+p)}dx\\ \label{PS3}&=&\Phi^{(d)}_{\kappa-m_d}(s^{-1})(\det s)^{-p-m_d}\\ \label{PS4}&=&\Phi^{(d)}_{\kappa}(s^{-1})(\det s)^{-p},
\end{eqnarray}
where  \eqref{PS2} and  \eqref{PS4} come from  \eqref{PS1},  and  \eqref{PS3}  from  \eqref{BZ},  where  $p$ is replaced by  $p+m_d.$

 From \eqref{EXPTR} we know that, for $2p\geq d-1$, the Laplace transform of $m(2p,d,d)$ is
\begin{equation}\label{LTM1}\int_{\mathcal{P}_{d}}e^{-\tr(sx)}m(2p,d,d)(dx)=(\det s)^{-p}\sum_{\kappa\in \mathcal{E}_d}\frac{C^{(d)}_{\kappa}(s^{-1})}{|\kappa|!}\;.\end{equation}
Observe that the Laplace transform of $f_d(t)\textbf{1}_{\mathcal{P}_d}(t)dt$ as defined by \eqref{ABS6} is easily deduced from \eqref{BZ} and is  equal to
\begin{equation}
\label{LTg}
\int_{\mathcal{P}_{d}}e^{-\tr(st)}f_d(t)dt=(\det s)^{-(d-1)/{2}}\left\{\sum_{\kappa\in \mathcal{E'}_{d}}\frac{C^{(d)}_{\kappa}(s^{-1})}{|\kappa|!}\right\}.
\end{equation}
In \eqref{LTM1}, take $2p=d-1$. Using the Laplace transform  \eqref{LTg},  we now want to prove that the Laplace transform of $r(dt)$ is 
\begin{equation}\label{LTM4}\int_{\mathcal{P}_{d}}e^{-\tr(st)}r(dt)=(\det s)^{-(d-1)/2}\left\{\sum_{\kappa\in\mathcal{E}_{d}\setminus \mathcal{E'}_{d}}\frac{C^{(d)}_{\kappa}(s^{-1})}{|\kappa|!}\right\}=(\det s)^{-(d-1)/2}\left\{\sum_{\kappa\in\mathcal{E}_{d-1}}\frac{C^{(d)}_{(\kappa,0)}(s^{-1})}{|\kappa|!}\right\}\;.\end{equation}
To prove  \eqref{LTM4} we start from the  definition of $r(dt)$. Observe first that, for $2p>d-1$, \eqref{BZ} and \eqref{LTM1} imply  
\begin{equation}\label{LTM2}m(2p,d,d)(dx)=(\det x)^{p-(d+1)/2}\left\{\sum_{\kappa\in \mathcal{E}_d}\frac{C^{(d)}_{\kappa}(x)}{|\kappa|! \, \Gamma_d(\kappa+p)}\right\}\textbf{1}_{\mathcal{P}_{d}}(x)dx.\end{equation}
In particular, in \eqref{LTM2} let us replace $d$ by $d-1$ and let $2p=d-1.$ We obtain 
\begin{equation}\label{LTM3}
(\det x)^{1/2}m(d-1,d-1,d-1)(dx)=\left[\sum_{\kappa\in \mathcal{E}_{d-1}}\frac{C^{(d-1)}_{\kappa}(x)}{|\kappa|! \, \Gamma_{d-1}\{\kappa+(d-1)/2\}}\right] \textbf{1}_{\mathcal{P}_{d-1}}(x)dx.\end{equation}
We can now write 
 \begin{eqnarray}
 \int_{\overline{\mathcal{P}}_{d}}e^{-\tr(st)}r(dt)&=& \label{L-1}\frac{\pi ^{1/2}}{\Gamma(d/2)}\int_{\mathbb{O}(d)}
 \left[ \int_{\mathcal{P}_{d-1}}e^{-\tr(su\widetilde{x}u^\top)}\sum_{\kappa\in \mathcal{E}_{d-1}}\frac{C^{(d-1)}_{\kappa}(x)}{|\kappa|! \, \Gamma_{d-1}\{ \kappa+(d-1)/2\} }dx\right] du
\\&=&\nonumber \frac{\pi ^{1/2}}{\Gamma(d/2)}\sum_{\kappa\in \mathcal{E}_{d-1}}\frac{C^{(d-1)}_{\kappa}(I_{d-1})}{|\kappa|!}\int_{\mathbb{O}(d)}\left[ \int_{\mathcal{P}_{d-1}}e^{-\tr([u^\top su]_1x)}\frac{\Phi^{(d-1)}_{\kappa}(x)}{\Gamma_{d-1}\{\kappa+(d-1)/2\}}dx\right] du.
\end{eqnarray}
Equality (\ref {L-1}) comes from (\ref {LTM3}) and the definition of $r$. 
Now we compute the last double integral as follows \begin{eqnarray}
\int_{\mathbb{O}(d)}\left[ \int_{\mathcal{P}_{d-1}}e^{-\tr([u^\top su]_1x)}\frac{\Phi^{(d-1)}_{\kappa}(x)}{\Gamma_{d-1}\{ \kappa+(d-1)/2\} }dx\right] du
&=&\label{L0}\int_{\mathbb{O}(d)}(\det[usu^\top]^{-1}_1)^{ {d}/{2}}\Phi^{(d-1)}_{\kappa}([usu^\top]^{-1}_1)du
\\\label{L1}&=&\int_{\mathbb{O}(d)}\Phi^{(d-1)}_{\kappa+{d}/{2}}([usu^\top]^{-1}_1)du
\\&=&\label{L2}\int_{\mathbb{O}(d)}\Phi^{(d-1)}_{-m_{d-1}-{d}/{2},\ldots,-m_{1}-{d}/{2}}([usu^\top]_1)du
\\&=&\label{L3}\Phi^{(d)}_{-m_{d-1}- {d}/{2},\ldots,-m_{1}- {d}/{2},0}(s)
\\&=&\label{L4}\Phi^{(d)}_{0,m_{1}+ {d}/{2},\ldots,m_{d-1}+ {d}/{2}}(s^{-1})
\\&=&\label{L5}\Phi^{(d)}_{- {d}/{2},m_{1},\ldots,m_{d-1}}(s^{-1})(\det s^{-1})^{ {d}/{2}}
\\&=&\label{L6}\Phi^{(d)}_{m_{1},\ldots,m_{d-1},0}(s^{-1})(\det s^{-1})^{(d-1)/2}.
\end{eqnarray} 

In the equations above, equality (\ref {L0}) follows from (\ref {BZ}) by replacing $(d,p)$ by $(d-1,d/2).$ Equalities (\ref {L1}) and \eqref{L5} come from Lemma~3. In the identities following (\ref {L2}) we have replaced $\kappa$ by $(m_1,\ldots,m_{d-1})$ for clarity. Formulas (\ref {L2}) and  \eqref{L4} come from Lemma~2, and \eqref{L3} comes from Lemma~1. The proof of  (\ref {L6}) is more involved and is a consequence of formula (iii) in~Theorem XIV 3.1 of FK, where we replace $(d,r,\lambda,\mu)$ respectively by $1,d$ and 
\begin{eqnarray*}\lambda &=&\Big(m_1+\frac{d-1}{4},m_2+\frac{d-3}{4},\ldots,m_{d-1}-\frac{d-3}{4},-\frac{d-1}{4}\Big),\\ \mu &=&\Big(-\frac{d-1}{4},m_1+\frac{d-1}{4},m_2+\frac{d-3}{4},\ldots,m_{d-1}-\frac{d-3}{4}\Big)\;.\end{eqnarray*} The fact that $\mu$ is a permutation of $\lambda$ and the reference above imply (\ref {L6}). 
Now we observe that  
\begin{equation}
\label{L7}
\frac{\pi ^{1/2}}{\Gamma(d/2)} \, C_{\kappa}(I_{d-1})=C_{\kappa,0}(I_{d}),
\end{equation}  
implied by formula \eqref{CST2}.  Finally we gather \eqref{L-1}, \eqref{L6} and \eqref{L7} to obtain 
$$
\int_{\overline{\mathcal{P}}_{d}}e^{-\tr(st)}r(dt)=(\det s)^{-(d-1)/2}\sum_{\kappa\in\mathcal{E}_{d-1}}\frac{C_{(\kappa,0)}(s^{-1})}{|\kappa|!} \, ,
$$
which proves \eqref{LTM4} and  Proposition~9 itself. \hfill $\square$

\subsection{Example: $m(2,3,3)$} 
\label{sec:4.4}

To illustrate Proposition~9 we consider the function $f_3(t)$ defined on $\mathcal{P}_3$ by \eqref{ABS6}. More specifically we have 
$$
f_3(t)=\frac{1}{\det t}\sum_{m_1\geq m_2\geq m_3>0}\frac{C^{(3)}_{m_1,m_2,m_3}(t)}{(m_1+m_2+m_3)! \, m_1! \,  \Gamma(m_2+{1}/{2}) \, (m_3-1)!}.
$$
We also consider the measure $ m(2,2,2)(dt)$ on $\mathcal{P}_2$ parameterized by $ (a,b,c)\mapsto t=\varphi(a,b,c)$ as in \eqref{VP}. From  \eqref{LTM3}, we have 
\begin{eqnarray*}m(2,2,2)(da,db,dc)&=&\frac{1}{(\det t)^{1/2}}\sum_{m_1\geq m_2\geq 0}\frac{C^{(2)}_{m_1m_2}(t)}{(m_1+m_2)! \, m_1! \, \Gamma(m_2+{1}/{2})}\\&=&
\frac{1}{\sqrt{\pi}}\sum_{k=0}^{\infty}\sum_{n=0}^{\infty}\frac{(a^2-b^2-c^2)^{n+(k-1)/{2}}2^{2n}}{({3}/{2}+k)_n(n+k)! \,  k! \, (2n)!} \, P_k\left(\frac{a}{\sqrt{a^2-b^2-c^2}}\right).
\end{eqnarray*} 

This last formula is obtained by using the calculations done in Section 4.2 for $\Phi^{(2)}_{m_1,m_2}$ in \eqref{LP2}, for $C^{(2)}_{m_1,m_2}(I_2)$ 
 in \eqref{CON2} and the change of indexes $(m_1-m_2, m_1)=(k,n).$
 Finally the singular measure $r(dt)$ concentrated on the set of matrices of rank 2 in the cone $\overline{\mathcal{P}_3}$ of  positive semidefinite matrices of order 3 is constructed as follows. One considers the product of  $m(2,2,2)(da,db,dc)$ by  the uniform probability measure $du$ on the orthogonal group $\mathbb{O}(3)$. The   measure $r(dt)$ is the image of this product measure by the following map: 
$$
(a,b,c,u)\mapsto t=u\left[\begin{array}{ccc}a+b&c&0\\c&a-b&0\\0&0&0\end{array}\right]u^\top.
$$
 Proposition~9 says that the measure $m(2,3,3)(dt)$  on the set $\overline{\mathcal{P}}_3$ of semipositive definite matrices of order 3 defined by the Laplace transform 
$ (\det s)^{-1}\exp \mathrm{trace} \ s^{-1}$ is the sum $r(dt)+f_3(t) \mathbf{1}_{\mathcal{P}_3}(t)dt.$

\section{Convolution lemmas in the cone $\overline{\mathcal{P}_d}$}  
\label{sec:5}

   Let $\mathcal{G}_k$ the set of linear subspaces $G$ of dimension $k$ of a Euclidean space $E$ of dimension $d$. Let us endow $\mathcal{G}_k$
with the uniform distribution, i.e., the unique probability  on $\mathcal{G}_k$ such that $G\sim uG$ for all $u\in \mathbb{O}(d).$  Lemma~4 below describes an intuitively obvious fact. For the sake of completeness, we offer a proof, although other ones may already exist in the literature. 

%%% Former Lemma 5.1
\bigskip \noindent \textbf{Lemma~4.}  
{\it Let $F$ be a fixed linear subspace of dimension $n$ of the Euclidean space $E$ of dimension $d.$ If the random linear subspace $G$ of $E$ has the uniform distribution on  $\mathcal{G}_k$ and if $k\leq d-n$, then $\Pr(G\cap F\neq \{0\})=0$.} 

\bigskip \noindent 
\textit{Proof.} It is enough to prove the lemma for $E=\mathbb{R}^d$, $F=\{0\}\times \mathbb{R}^{d-n}$ and $k=d-n$.  Let $Z_1,\ldots,Z_n$ be iid random variables in $\mathbb{R}^d$ following the standard Gaussian distribution $\mathcal{N}(0,I_d)$. Let $G$ be the random linear subspace of $E$ generated by the vectors $Z_1,\ldots,Z_n.$ Since for all $u\in \mathbb{O}(d)$ we have $(uZ_1,\ldots,uZ_n)\sim (Z_1,\ldots,Z_n)$, clearly $G\sim uG$ and $G$ follows the uniform distribution. Introduce the matrix 
$$
M=[Z_1,\ldots,Z_n]=(Z_{ij})_{1\leq i\leq d,\ 1\leq j\leq n},
$$
 whose columns are the vectors $Z_1,\ldots,Z_n.$ Then $x_1Z_1+\cdots+x_nZ_n=MX$, where $X=(x_1,\ldots,x_n)^\top$.
Now $G\cap F\neq \{0\}$ implies that there exists a non-zero $X$ such that the first $n$  elements of $MX$ are zero. In other terms, considering the square matrix $M_1$ of order $n$ defined by $M_1=(Z_{ij})_{1\leq i,j\leq n}$, we have that $G\cap F\neq \{0\}$ implies that  there exists a non-zero $X$ such that $M_1X=0$. This happens if and only if $\det M_1=0.$ Since the $n^2$ entries of the matrix $M_1$ are independent $\mathcal{N}(0,1)$ variables, the event $\det M_1=0$ has probability zero and  the lemma is proved. \hfill $\square$

\bigskip
In the sequel we will denote by $S_b$ the set of $x\in \overline{\mathcal{P}_d}$ with  $b=\rank x=0,\ldots,d.$ Of course $S_d=\mathcal{P}_d.$

%%%Former Lemma 5.2
\bigskip \noindent 
\textbf{Lemma~5.} 
{\it Let $Y$ be a random variable in $S_b$ and assume that $uYu^\top\sim Y$  for all $u$ in the orthogonal group $\mathbb{O}(d).$ Let  $x_0\in S_a.$ Then $x_0+Y$ is concentrated on $S_{a+b}$ if $a+b\leq d$ or on $S_d=\mathcal{P}_d$ if $a+b\geq d.$  Furthermore if $x_0\in S_c,$  if $x_0+Y$ is concentrated on $S_{a+b}$ and if $a+b<d$, then $c=a$.} 

\bigskip \noindent \textbf{Remark 3.} If $a+b=d$ and if $x_0+Y$ is concentrated on $S_{a+b}=S_d,$ $x_0$ could be on any $S_c$ with $a\leq c\leq k.$ 

\bigskip \noindent \textit{Proof.} Apply Lemma~4 to $F=x_0\mathbb{R}^d$ and to $G=Y\mathbb{R}^d.$ Then almost surely we have  $\dim (F+G)=a+b$ if $a+b\leq d$. Furthermore we  always have $\rank(x_0+Y)\leq a+b.$ 

To see that $\rank(x_0+Y)= a+b$ almost surely, let us suppose that $(x_0+Y)\mathbb{R}^d\neq E=F+G.$ Let $x_0'$ and $Y'$ be the restrictions of the endomorphisms $x_0$ and $Y$ to the linear space $E.$ Since $x_0$ and $Y$ are symmetric, this implies that $x_0'E=F$ and $Y'E=G.$ Since  $(x'_0+Y')E\neq E$, there exists $v\in E\setminus\{0\}$ which is orthogonal to $(x'_0+Y')E$ and thus $(x'_0+Y')v=0.$ Since $x'_0v\in F$ and $Y'v\in G$ and since $F\cap G=\{0\}$, this implies that $x'_0v=Y'v=0$, and $v$ is in $\mathrm{Ker}(x'_0)\cap\mathrm{Ker}(Y').$ But since we have almost surely $F\oplus G=E$ (a direct sum, not necessarily an orthogonal one), we have also almost surely $\mathrm{Ker}(x'_0)\oplus \mathrm{Ker}(Y')=E$, which implies $\mathrm{Ker}(x'_0)\cap\mathrm{Ker}(Y')=\{0\}.$ Thus almost surely  $v=0$, which is a
contradiction. Finally, $(x_0+Y)\mathbb{R}^d=E=F+G$ and $\rank(x_0+Y)= a+b.$ If $a+b>d$ then $F$ contains a subspace of dimension $d-b$ and $\dim (F+G)=d.$ 

To conclude, suppose now that $x_0\in S_c,$ and that $x_0+Y$ is concentrated on $S_{a+b}$ with $a+b<k.$  If $c+b<k$ then, by the first part of the lemma,  $x_0+Y$ is concentrated on $S_{c+b}$. But $S_{a+b}=S_{c+b}$ implies $c=a.$  If $c+b\geq k$ then, by the first part of the lemma again, 
$x_0+Y$ is concentrated on $S_{d}$. This is impossible since $S_{a+b}\neq S_{d}.$ \hfill $\square$

\bigskip
\noindent 
%%% Former Lemma 5.3
\textbf{Lemma~6.} 
{\it Let $\mu$ and $\nu$ be positive measures on $\overline{\mathcal{P}_d}$ such that $\nu$ is concentrated on $S_b$ and $\nu$ is invariant by  the transformations $x\mapsto uxu^\top,\;u\in \mathbb{O}(d)$. Let $a \in \{0,\ldots,d-2\}$ be such that $a+b<d.$  Then  $\mu*\nu$ is concentrated on $S_{a+b}$ if and only if $\mu$ is concentrated on $S_a.$  Furthermore $\mu*\nu$ is concentrated on $S_{d}=\mathcal{P}_d$
if $\mu$ is concentrated on $S_{d-b}$.} 

\bigskip
\noindent \textit{Proof.} $\Rightarrow$: For $y_0\in S_b$ consider the distribution $K_{y_0}(dy)$ on $S_b$ of the random variable $Uy_0U^\top$, where $U$ is uniformly distributed on the orthogonal group $\mathbb{O}(k).$ Let $D_b$ the set of diagonal elements $y_0$ of $S_b$ of the form $y_0=\mathrm{diag}(\lambda_1,\ldots,\lambda_b,0,\ldots,0)$ such that $\lambda_1\geq \cdots\geq \lambda_b>0.$ Then there exists a unique positive measure $\nu_0$ on $D_b$
such that the following disintegration holds:
$$
\nu(dy)=\int_{D_b}\nu_0(dy_0)K_{y_0}(dy).
$$ 
It follows that
\begin{eqnarray*}(\mu*\nu)(dx)=\int_{S_b}\nu(dy)\mu(dx-y)=\int_{S_b}\mu(dx-y)\int_{D_b}\nu_0(dy_0)K_{y_0}(dy) =\int_{D_b}\nu_0(dy_0)\int_{S_b}\mu(dx-y)K_{y_0}(dy).
\end{eqnarray*}
Therefore the measure $\mu*K_{y_0}$ is concentrated on $S_{a+b}$ for  $ \nu_0$ almost all $y_0\in D_b$. From Lemma~5, this implies that $\mu$ is concentrated on $S_a$.

${\Leftarrow}:$ If $\mu$ is concentrated on $S_a$ with $a+b\leq d$, it is an easy consequence of Lemma~5 that $\mu*\nu$ is concentrated on $S_{a+b}.$ \hfill $\square$

%%% Former Lemma 5.4
\bigskip \noindent \textbf{Lemma~7.} 
{\it Let $a,b\in \{1,\ldots,d-1\}$ such that $a+b<d.$ If $m(a,a+b,d)$  exists, it is concentrated on $S_a$.}

\bigskip 
\noindent 
\textit{Proof.} From the  Laplace transforms  of the measures we know that $m(a,a+b,d)*m(b,0,d)=m(a+b,a+b,d)$. From Proposition~2.2, we know that $ m(a+b,a+b,d)$  is concentrated on $S_{a+b}$ since this is the case for the singular non-central  Wishart $\mathcal{N}(n,I(n,d),I_d)$ with $n=a+b.$ Since the Laplace transform of $m(b,0,d)$ is $(\det s)^{-b/2}$, we know that $m(b,0,d)$ is invariant by the transformations $x\mapsto uxu^\top$ for any $u\in \mathbb{O}(d).$ By Lemma~6 we deduce that $m(a,a+b,d)$ is concentrated on $S_{a}$ if it exists. \hfill $\square$

\section{$m(d-2,d-1,d)$ and $m(d-2,d,d)$ do not exist for $d\geq 3$}
\label{sec:6}

In this section we prove Propositions~4 and 5. 

\bigskip \noindent 
\textit{Proof of Proposition~4.}  Suppose that $m(d-2,d,d)$  exists. Then $m(d-2,d,d)*m(2,0,d)=m(d,d,d).$ From  \eqref{BZ} and \eqref{LTM1}, the measure $m(d,d,d)$ has a $C^{\infty}$  density $g.$  As a consequence
$$
\int_{\overline{\mathcal{P}_d}}e^{-\tr( sx)}m(d-2,d,d)(dx)=\det s \int_{\overline{\mathcal{P}_d}}e^{-\tr( sx)}g(x)dx=(-1)^n\int_{\overline{\mathcal{P}_d}}e^{-\tr (sx)}\det( {\partial}/{\partial x})g(x)dx.
$$ 
This implies that  $m(d-2,d,d)(dx)=(-1)^n\det({\partial}/{\partial x})g(x)dx$ has a density. However,  
since $m(d-2,d,d)*m(1,0,d)=m(d-1,d,d)$ this would imply that $m(d-1,d,d)$ is absolutely continuous, which contradicts Proposition~9. \hfill $\square$
 
\bigskip \noindent 
\textit{Proof of Proposition~5.} Suppose that $m(d-2,d-1,d)$  exists.
By Lemma~7, the measure $m(d-2,d-1,d)$ is concentrated on $S_{d-2}$. Therefore there exists a positive measure $m(dy)=m(dy_1,\ldots,dy_{d-2})$ on $\mathbb{R}^{d(d-2)}=\mathbb{R}^d\times\cdots \times \mathbb{R}^d $  such that for all $s\in \mathcal{P}_d$ we have
$$
\int_{\mathbb{R}^{d(d-2)}}e^{-(y_1^\top sy_1+\cdots+ y_{d-2}^\top sy_{d-2})}m(dy)=\frac{1}{\det s^{(d-2)/2}}e^{\tr \{s^{-1}I(d-1,d)\}}.
$$
We write  the elements $y=(y_1,\ldots,y_{d-2})$ more conveniently with the help of the transposed  matrix $y^\top=(y_{i,j})$ with $d-2$ rows $y_1^\top,\ldots,y_{d-2}^\top$ and $d$ columns $c_1,\ldots,c_d$
$$
y^\top=\left[\begin{array}{c}y_1^\top\\\vdots\\y_{d-2}^\top\end{array}\right]
=\left[\begin{array}{ccc}y_{1,1}&\dots&y_{1,d}\\\vdots&\vdots&\vdots\\y_{d-2,1}&\dots&y_{d-2,d}\end{array}\right]=\left[c_1,\dots,c_{d}\right].$$
With this notation introduce the Gram matrix $G(c)=G(c_1,\ldots,c_d)=(\<c_j,c_k\>)_{1\leq j,k\leq d}$ and denote by $m(dc)$ what we denoted by $m(dy)$ before. We get
\begin{equation}\label{SG}\int_{\mathbb{R}^{d(d-2)}}e^{-\tr \{sG(c)\}}m(dc)=\frac{1}{\det s^{(d-2)/2}}e^{\tr \{s^{-1}I(d-1,d)\}}.\end{equation}
Equality \eqref{SG} means that $m(d-2,d-1,d)(dx)$ is the image of $m(dc)$ by $c\mapsto x=G(c).$ 

Now in \eqref{SG} we choose $s=\mathrm{diag}(1,s_1)$ where $s_1$ is a symmetric positive definite matrix of order $d-1.$ 
We also desintegrate $m(dc)$ by introducing a probability kernel $K(c_2,\ldots,c_d;dc_1)$ and a positive measure $m_1(dc_2,\ldots,dc_d)$ such that 
$$
e^{-\|c_1\|^2}m(dc_1,dc_2,\ldots,dc_d)=m_1(dc_2,\ldots,dc_d)K(c_2,\ldots,c_d;dc_1).
$$ 
With these notations we can write
\begin{eqnarray*}
\frac{1}{\det s_1^{(d-2)/2}}e^{\tr (s_1^{-1})}& =& \int_{\mathbb{R}^{d(d-2)}}e^{-\tr \{sG(c)\}}m(dc) =\int_{\mathbb{R}^{d(d-2)}}e^{-\|c_1\|^2}e^{-\tr \{s_1G(c_2,\ldots,c_d)\}}m(dc_1,dc_2,\ldots,dc_d)\\
&=&\int_{\mathbb{R}^{(d-1)(d-2)}}e^{-\tr \{s_1G(c_2,\ldots,c_d)\}}\left\{\int_{\mathbb{R}^{d-2}}K(c_2,\ldots,c_d;dc_1)\right\} m_1(dc_2,\ldots,dc_d)\\
&=&\int_{\mathbb{R}^{(d-1)(d-2)}}e^{-\tr \{s_1G(c_2,\ldots,c_d)\}}m_1(dc_2,\ldots,dc_d)
\end{eqnarray*}
since $K$ is a probability kernel. 
The last equality says that the image of $m_1(dc_2,\ldots,dc_d)$ by the map $(c_2,\ldots,c_d)\mapsto x=G(c_2,\ldots,c_d)$ is nothing but $m(d-2,d-1,d-1)(dx).$ Denote $G_2=G(c_2,\ldots,c_d)$ for simplicity. Since $c_2,\ldots,c_d$ are vectors of a Euclidean space of dimension $d-2$, the rank of $ G_2$ is  less than or equal to  $d-2$. To prove this elementary fact of linear algebra we use   $G_2\in\overline{\mathcal{P}_{d-1}}.$ This implies that if $x=(x_2,\ldots,x_d)^\top$, then $G_2x=0$ if and only if $x^\top G_2x=0.$ Since $x^\top G_2x=\| x_2c_2 + \cdots + x_dc_d \|^2$ the linear space of $x\in \mathbb{R}^{d-1}$ such that $ x_2c_2 + \cdots + x_dc_d=0$ has at least dimension 1, the kernel of the endomorphism of  $\mathbb{R}^{d-1}$ with matrix $G_2$ has at least dimension 1 and its image has at most dimension $d-2.$ 
This contradicts Proposition~9 which says  that $m(d-2,d-1,d-1)$ has an absolutely continuous part and therefore charges matrices with rank $d-1$. \hfill $\square$

\section{Bibliographical comments and acknowledgments}
\label{sec:7}

The question of the existence of the non-central Wishart distribution was addressed in \cite{LM2008} where we claimed, in Proposition 2.3, that such a distribution exists if and only if $p$ is in $\Lambda_d$ 
without any restriction on $\Sigma\in \mathcal{P}_d$ or  $w\in \overline{\mathcal{P}_d}.$  However the proof of the `if' part was not given in \cite{LM2008}  since  we  considered it obvious that if $p$ was in the part $\{{1}/{2}, {2}/{2},\ldots, (d-1)/{2}\}$ of $\Lambda_d$, then $\mathcal{NCW}(2p,w,\Sigma)$ did exist without restriction on the rank of $w$. This gap in \cite{LM2008}  was kindly pointed out to us, in a private communication, by E. Mayerhofer who later showed in \cite{Mayerhofer}  that the statement was not only unproven, but false. More specifically, Mayerhofer \cite{Mayerhofer} showed that if $\mathcal{NCW}(2p,w,\Sigma)$ exists, if $d\geq 3$ and if $n=2p$ is in $\{1,2,\ldots,d-2\}$, then $\rank w\leq n+1$. We reproved it in a different form in Proposition~2.4 of the present paper.  Mayerhofer \cite{Mayerhofer}  used a stochastic process valued in the set of symmetric matrices in order to prove this statement.  Finally, he conjectured in  \cite{Mayerhofer} that $\rank w\leq n$  holds, and Proposition~2.5  shows that the conjecture is true. 

Our present  paper corrects the  mistake in \cite{LM2008} by giving a  necessary and sufficient condition for the existence of $\mathcal{NCW}(2p,w,\Sigma)$ through Proposition 2.6 and thus gives a proof of the aforementioned conjecture by  Mayerhofer. We  follow the lines of the arXiv paper \cite{LM2011}, with a correction to the proof of Proposition 2.4 and a complement to Proposition 2.1.  Recently, Graczyk, Ma\l ecki and Mayerhofer in \cite{GMM} have also given  a proof of the same result, using  the same techniques as  \cite{Mayerhofer}. Clearly, the methods of the present paper are based on linear algebra and are of a different nature. 

Because of the gap in the `if' part of Proposition 2.3 of  \cite{LM2008}, it happened that the proof of the `only if ' part became also incomplete. This fact has been pointed out by Piotr Graczyk in a private communication. With the help of Mauro Piccioni, we have been able to design the simple proof of the second part of Proposition 2.1 of the present paper. We have taken the idea of Proposition 2.3 from  \cite{Mayerhofer}. Jacques  Faraut helped  us with Lemma 4.2. 

The `only if' part of Proposition 2.6 has a long history: for $w=0$, an unconvincing proof appears in \cite{Olkin} and is commented upon in \cite{Casalis}.  Later on, this result was conjectured by Eaton \cite{Eaton}, who was unaware of \cite{Gindikin}. Gindikin's proof is described in FK Theorem VII 3.1, where it is explained why \eqref{GYN} is also called the Wallach set. Shanbhag in \cite{Shanbhag} gave an elementary and elegant proof, which is now the classical one. Peddada and Richards \cite{Peddada} gave the proof described in Proposition 4.1 of the present paper. They also  proved  that $p\in \Lambda$ when $w$ has rank one. 

The proper mathematical framework for this paper is that  of  Euclidean Jordan algebras rather than the linear spaces of real symmetric matrices. But working in  that general framework
  might have obscured our statements without adding any insight: the extension of our results to Euclidean Jordan algebras is straightforward.  However, past experience comparing the real and the complex case \cite{GLM2003,GLM2005} makes us aware  that replacing real symmetric matrices by Hermitian ones could  lead to more  elegant formulas than in the real case.  
To conclude, our deepest thanks go to Eberhard Mayerhofer, Jacques Faraut, Piotr Graczyk, Mauro Piccioni and the two referees.  

\bigskip
\noindent
\textbf{Acknowledgments}. H. Massam gratefully acknowledges support from the Natural Sciences and Engineering Research Council of Canada through Discovery Grant No A8947.

\end{document}